\documentclass[12pt, reqno]{amsart}

\usepackage{lmodern}
\usepackage{mathrsfs}
\usepackage{physics}
\usepackage{mathtools}
\usepackage{amsmath,amsfonts,latexsym,amssymb,mathrsfs}
\usepackage{comment,url,soul,microtype,color}
\usepackage{tikz, tikz-cd, graphicx, caption, float}
\usepackage[backref=page, bookmarks=false, 
colorlinks=true,
linkcolor=blue,
citecolor=blue,
filecolor=blue,
menucolor=blue,
urlcolor=blue]{hyperref}

\let\oldtocsection=\tocsection
\let\oldtocsubsection=\tocsubsection
\let\oldtocsubsubsection=\tocsubsubsection

\renewcommand{\tocsection}[2]{\hspace{0em}\oldtocsection{#1}{#2}}
\renewcommand{\tocsubsection}[2]{\hspace{3em}\oldtocsubsection{#1}{#2}}
\renewcommand{\tocsubsubsection}[2]{\hspace{6em}\oldtocsubsubsection{#1}{#2}}

\newtheorem{lemma}{Lemma}
\newtheorem{proposition}[lemma]{Proposition}
\newtheorem{definition/proposition}[lemma]{Definition/Proposition}
\newtheorem{corollary}[lemma]{Corollary}
\newtheorem{theorem}[lemma]{Theorem}

\theoremstyle{definition}
\newtheorem{definition}[lemma]{Definition}
\newtheorem{remark}[lemma]{Remark}
\newtheorem{observation}[lemma]{Observation}

\newtheorem{example}[lemma]{Example}
\newtheorem{conjecture}[lemma]{Conjecture}

\newtheorem*{example*}{Example}
\newtheorem{theorema}{Theorem}

\numberwithin{equation}{section}

\newcommand{\plusm}{\partial_{+}X} 
\newcommand{\minm}{\partial_{-}X}
\newcommand{\dist}{\mathrm{dist}}

\newcommand{\mts}{M\times S^1}

\newcommand{\dvols}{\,dvol_{\Sigma}}

\newcommand{\RR}{\mathbf{R}}
\newcommand{\ZZ}{\mathbf{Z}}
\newcommand{\DD}{\mathbf{D}}
\newcommand{\cpt}{\mathbf{CP}^3}
\newcommand{\cptwo}{\mathbf{CP}^2}
\newcommand{\cptwobar}{\overline{\mathbf{CP}}^2}

\DeclareMathOperator{\wid}{width}
\DeclareMathOperator{\Wh}{Wh}

\begin{document}

\title[Positive scalar curvature and exotic structures]{Positive scalar curvature and exotic structures on simply connected four manifolds}

\author{Aditya Kumar}
\address{Department of Mathematics, Johns Hopkins University. 3400 N. Charles Street, Baltimore, MD 21218, USA}
\email{akumar65@jhu.edu}

\author{Balarka Sen}
\address{School of Mathematics, Tata Institute of Fundamental Research. 1, Homi Bhabha Road, Mumbai-400005, India}
\email{balarka2000@gmail.com, balarka@math.tifr.res.in}



\begin{abstract}
We address Gromov's band width inequality and Rosenberg's $S^1$-stability conjecture for simply connected smooth four manifolds. Both results are known to be false in dimension $4$ due to counterexamples based on Seiberg-Witten invariants. Nevertheless we show that both of these results hold upon considering simply connected smooth four manifolds up to homeomorphism. We also obtain a related result for non-simply connected smooth four manifolds.
\end{abstract}

\maketitle
\tableofcontents

\section{Introduction}

In this article, we are concerned with simply connected $4$-manifolds admitting a metric of positive scalar curvature (PSC). The class of closed manifolds admitting PSC metrics is completely understood in dimension $2$ and $3$. Therefore, all the problems discussed in this paper are known when $n=2,3$.

We will study two problems for manifolds with positive scalar curvature, namely, Rosenberg's $S^1$-stability conjecture \cite[Conjecture 1.24]{ros07} and Gromov's band width inequality conjecture \cite[11.12, Conjecture C]{gromovmetric}. 

\begin{conjecture}[$S^1$-stability]\label{conj:S1stability}
    Let $M^n$ be a closed manifold. Then $M$ admits a PSC metric if and only if $\mts$ admits a PSC metric.
\end{conjecture}
\begin{conjecture}[Width inequality]\label{conj:width}
    Let $M^n$ be a closed manifold that does not admit a PSC metric. Consider the band $X^{n+1} = M^n \times [-1,+1]$ equipped with a metric $g$, such that the scalar curvature satisfies the lower bound $R_g \geq n(n+1)\kappa^2$. Then, 
    \begin{equation*}
                    \wid(X,g) := \dist_g(M \times \{-1\}, M \times \{1\})    \leq \frac{2\pi}{n+1}\frac{1}{\kappa}
    \end{equation*}
\end{conjecture}

Note that the width inequality implies the $S^1$-stability conjecture by passing to the infinite cyclic cover $M\times \RR$. Therefore, a counterexample to the $S^1$-stability conjecture is also a counterexample to the width inequality conjecture.

In the \emph{Four Lectures} \cite[Section 5]{gromovfour}, Gromov outlined a proof for Conjecture \ref{conj:width}, which was subsequently completed by R{\"a}de in \cite{rade23} for dimension $n=5,6$. As noted previously, it was already known in dimension $n=2,3$ for other reasons.

Note that in dimensions $2$ and $3$ the only simply connected manifolds are $S^2$ and $S^3$ respectively, both of which have positive scalar curvature with respect to the standard round metric. In dimensions $5$ and above there is the classification result of Gromov-Lawson-Stolz for simply connected manifolds admitting a metric of positive scalar curvature.

\begin{theorem}[Gromov-Lawson, Stolz \cite{glsc,stolz}] \label{theorem:gls}
    For $n\geq 5$, let $M^n$ be a closed simply connected manifold. If $M$ is not spin, then it admits a metric of positive scalar curvature \cite{glsc}. If $M$ is spin, then it admits a metric of positive scalar curvature if and only if $\hat{\alpha}(M)=0$ \cite{stolz}.
\end{theorem}
\begin{remark}
When $n=4$ as in this paper, $\hat{\alpha}(M)$ is just the Lichnerowicz genus $\hat{A}(M)$ that is given by $-\sigma(M)/8$, where $\sigma(M)=b_2^+-b_2^-$ is the signature. $\hat{\alpha}$ was introduced by Milnor \cite{milnorahat} as a generalisation of $\hat{A}$ when $n \not \equiv 0 \pmod{4}$.  
\end{remark}

\subsection{Counterexamples to $S^1$-stability in dimension $4$} When $n=4$, these conjectures, and several other related results are known to be false due to an obstruction to admitting a PSC metric coming from non-vanishing of the Seiberg-Witten invariant. The following counterexample to $S^1$-stability in dimension $4$ was provided in Rosenberg. We will discuss it in detail as understanding this example was a primary motivation behind this work.

\begin{example*}[Counterexample to $S^1$ stability in $n=4$]\cite[p.23, Remark 1.25]{ros07}
Let $M^4 = V_5$, where $V_5$ is the zero set of a degree $5$ homogeneous polynomial in $\cpt$:
  $$V_5= \{[z_0:z_1:z_2:z_3] \in \cpt \, | \,  z_0^5+z_1^5+z_2^5+z_3^5 = 0 \} .$$ 
$V_5$ has $b_2^+=9$ and $b_2^-=44$. Since it is an algebraic hypersurface in $\cpt$, it has an induced symplectic structure. This along with $b_2^+(V_5) > 1$ implies that $V_5$ has non-zero Seiberg-Witten invariant due to a theorem of Taubes \cite{taub}. Consequently, $V_5$ does not admit a PSC metric. On the other hand, $V_5$ is simply connected, therefore $\pi_1(V_5 \times S^1)= \ZZ$. Further, $V_5 \times S^1$ is not spin as $V_5$ has odd degree. Therefore, by a theorem of Stolz, $V_5 \times S^1$ admits a PSC metric, even though $V_5$ does not. 
\end{example*}

While this argument tells us that $V_5 \times S^1$ admits a PSC metric, the provenance of such a metric is not clear. In general if $M$ admits a PSC metric, then $g_M+dt^2$ is a PSC metric on $M\times S^1$. Our starting observation is that the situation for $V_5 \times S^1$ is not too different. More explicitly, one can obtain a PSC metric on $V_5 \times S^1$ from an \emph{exotic smooth structure} on $V_5$ that admits a PSC metric. 
  
\begin{example*}[Direct argument that $V_5 \times S^1$ is PSC]
    As $V_5$ is a smooth simply connected $4$-manifold that is not spin, it has an odd intersection form that is isomorphic to $Q_{V_5} \cong b_2^+\langle 1 \rangle \oplus b_2^- \langle -1 \rangle$ \cite[Theorem 1.2.21]{gsbook}. For $V_5$, one had $b_2^+=9$ and $b_2^-=44$, which gives us $Q_{V_5} \cong 9\langle 1 \rangle \oplus 44 \langle -1 \rangle$. This is also the intersection form of $9\cptwo \# 44\cptwobar$.  We saw that $V_5$ does not admit a PSC metric due to Taubes' theorem. Nevertheless, $9\cptwo \# 44\cptwobar$ admits a PSC metric, since a connect sum of PSC manifolds of dimension at least $3$ is still PSC by the Gromov-Lawson-Schoen-Yau surgery theorem \cite{glsc,sydescent} .
    
    Since $V_5$ and $9\cptwo \# 44\cptwobar$ are both simply connected with isomorphic intersection forms, by Wall's theorem \cite{wall} $V_5$ is $h$-cobordant to $9\cptwo \# 44\cptwobar$. Thus, $V_5  \times S^1$ is also $h$-cobordant to  $(9\cptwo \# 44\cptwobar) \times S^1$. This is also an $s$-cobordism, since $\pi_1(V_5 \times S^1) = \ZZ$, and the Whitehead group of $\ZZ$ is trivial, $\Wh(\ZZ) = 0$. As this is a $6$-dimensional $s$-cobordism, we can apply the $s$-cobordism theorem \cite[p. 15]{luck}. Therefore, $V_5 \times S^1$ is diffeomorphic to $(9\cptwo \# 44\cptwobar) \times S^1$. Pulling back the PSC metric on $(9\cptwo \# 44\cptwobar) \times S^1$ by this diffeomorphism, one gets a PSC metric on $V_5 \times S^1$.
\end{example*}
\begin{remark}
    In general the same argument shows the following: Let $M$ be a simply connected Kahler surface. Then for every integer $k> 0$, $(M \# k\cptwobar) \times S^1$ admits a PSC metric. On the other hand if we also have $b_2^+(M)>1$, then $M \# k\cptwobar$ does not admit a PSC metric. Indeed, since $b_2^+(M)>1$ and blowing up\footnote{The connect sum $M \# \cptwobar$ is called a blow-up of $M$.} does not change $b_2^+$, $M \# k\cptwobar$ has non-zero Seiberg-Witten invariant by Taubes' theorem. However, performing a blow-up makes the intersection form odd. Therefore, $M \# k\cptwobar$ has an exotic copy $m \cptwo \# n \cptwobar$ which admits a PSC metric, where $m=b_2^+(M), n=b_2^-(M)+k$. Therefore, since $M\#k\cptwobar$ is simply connected, $(M\#k\cptwobar) \times S^1$ admits a PSC metric by the $s$-cobordism argument, as in the above example. 
\end{remark}

These counterexamples lead us to speculate that the failure of $S^1$ stability and several other related results on manifolds with positive scalar curvature in dimension $4$ is \textit{only} due to exotic structures. That is to say, if one looks at smooth four manifolds up to homeomorphisms, then the existing results could be extended to dimension $n=4$. In this paper we take this point of view and establish Gromov's width inequality and Rosenberg's $S^1$-stability conjecture for simply connected smooth four manifolds. We will also observe that this framework automatically extends the Gromov-Lawson-Stolz classification of simply connected PSC manifolds to dimension $4$. 

\subsection{Statements of main results} 
\begin{definition}
    We will say that a manifold $M$ is PSC if it admits a positive scalar curvature metric. We will say that a $4$-manifold $M$ is \emph{PSC upto homeomorphism}, if there exists at least one smooth structure on $M$ that admits a positive scalar curvature metric.
\end{definition}

The main results of the article are the following.

\begin{theorema}[Width inequality in dimension $4$]\label{theoremB}
    Let $M^4$ be a closed simply connected smooth $4$-manifold that is not PSC upto homeomorphism. Let $g$ be a metric on the band $X = M^4 \times [-1,1]$ such that $R_g \geq 20\kappa^2$. Then,
$$\wid(X^5,g) := \dist_g(M \times \{-1\},M \times \{1\}) \leq \frac{2\pi}{5\kappa}$$
\end{theorema}

As a consequence of Theorem \ref{theoremB} we obtain $S^1$-stability:

\begin{theorema}[$S^1$-stability in dimension $4$]\label{theoremA}
    Let $M^4$ be a closed simply connected smooth $4$-manifold. Then, $M$ is PSC upto homeomorphism if and only if $M \times S^1$ is PSC.
\end{theorema}
\begin{remark}
    Note that in the forward direction the result of Theorem \ref{theoremA} is slightly stronger than the statement of Conjecture \ref{conj:S1stability} because we do not assume that $M$ is PSC, only that it is PSC upto homeomorphism. 
\end{remark}

We will prove Theorem \ref{theoremB} as a consequence of a more general width inequality that we can show for all $4$-manifolds. 

\begin{theorema}\label{theoremC}
    Let $M$ be a closed smooth $4$-manifold such that $M\#^k(S^2 \times S^2)$ does not admit a PSC metric for any $k\geq 0$. Let $g$ be a metric on the band $X = M^4 \times [-1,1]$ such that $R_g \geq 20\kappa^2$. Then,
$$\wid(X^5,g) := \dist_g(M \times \{-1\},M \times \{1\}) \leq \frac{2\pi}{5\kappa}$$
 \end{theorema}

Theorem \ref{theoremC} has the following corollary. 

\begin{corollary}
   Let $M^4$ be any closed smooth $4$-manifold. If $\mts$ is PSC then $M\#^k(S^2 \times S^2)$ is PSC for some $k\geq 0$.
\end{corollary}

The key to passing from Theorem \ref{theoremC} to Theorem \ref{theoremB} will be the following observation. It will be deduced from standard facts about simply connected $4$-manifolds and the fact that the Lichnerowicz genus $\hat{A}$ is a homotopy invariant in dimension $4$. This is because $\hat{A}(M^4) = -\sigma(M)/8$, where $\sigma(M)=b_2^+(M) -b_2^-(M)$ is the signature of $M$.

\begin{observation}[Gromov-Lawson-Stolz in dimension $4$] \label{theoremCD}
     Let $M^4$ be a closed simply connected smooth $4$-manifold. If $M$ is not spin, then it is PSC upto homeomorphism. If $M$ is spin, then it is PSC upto homeomorphism if and only if $\hat{A}(M)=0$.
\end{observation}

Theorem \ref{theoremB} combined with Observation \ref{theoremCD} gives the following corollary. 

\begin{corollary}
    $K3 \times [-1,1]$ satisfies the width inequality. 
\end{corollary}

We close with an example of a non-spin, non-simply connected smooth four manifold $M^4$ for which we can deduce using the above results that $M \times [-1,1]$ satisfies the width inequality. 

\begin{example}
    Take $M^4 = K3\#X$ with $X = (S^2 \times S^2)/\ZZ_2$, where $\ZZ_2$ acts on $S^2 \times S^2$ by $\overline{1} \cdot (x, y) = (-x, -y)$. Note that $\pi_1(M)=\ZZ_2$. We will show that $M$ is not spin, and that $M\times[-1,1]$ satisfies the width inequality. 
    
    We first observe $X$ is not spin. Indeed, consider the diagonal embedding $\mathrm{diag} : S^2 \hookrightarrow S^2 \times S^2$, given by $\mathrm{diag}(x) = (x, x)$. Since $\mathrm{diag}$ is $\ZZ_2$-equivariant, it descends to an embedding
    $$\overline{\mathrm{diag}} : \mathbf{RP}^2 \hookrightarrow X$$
    The normal bundle of $\mathrm{diag}$ is the tangent bundle $TS^2$. Therefore, the normal bundle of $\overline{\mathrm{diag}}$ is the tangent bundle $T\mathbf{RP}^2$. If $X$ was spin, we would have $w_2(TX) = 0$. However, 
    \begin{align*}
        w_2(\mathrm{diag}^*TX) 
        = w_2(T\mathbf{RP}^2 \oplus N\mathbf{RP}^2) 
        &= w_2(T\mathbf{RP}^2 \oplus T\mathbf{RP}^2) \\
        &= 2w_2(T\mathbf{RP}^2) + w_1(\mathbf{RP}^2)^2 \\
        &= x^2
    \end{align*}
    where $H^*(\mathbf{RP}^2; \ZZ_2) \cong \ZZ_2[x]/(x^3)$. Therefore, $w_2(\mathrm{diag}^*TX) \neq 0$. This shows $X$ is not spin.

    We now show the band inequality for $M \times [-1,1]$. Note that $\pi_1(M) = \ZZ_2$ and $M$ is not spin. Nevertheless, the universal cover $\widetilde{M} = \#^2 K3 \# (S^2 \times S^2)$ is spin. If $M\times [-1,1]$ does not satisfy the band inequality then by Theorem \ref{theoremC}, $M \#^k (S^2 \times S^2)$ is PSC for some $k \geq 0$. Therefore, the universal cover 
    $$\widetilde{M} \#^{2k} (S^2 \times S^2) = \#^2 K3 \#^{2k+1} (S^2 \times S^2),$$
    must also be PSC. A spin PSC manifold must have $\widehat{A}=0$ \cite{lic63}. But here, $\widehat{A} = -4$. Contradiction.
\end{example}

\subsection{Idea of proof} The proof of Theorem \ref{theoremB} will be accomplished in three steps, following the same general scheme as  \cite[Section 5]{gromovfour} and \cite{rade23}.

Suppose a band $M^n \times [-1,1]$ admits a positive scalar curvature metric. Then, if the band is long enough, i.e., if $\dist_g(M\times\{-1\},M\times\{1\})$ is above a certain threshold, then the band contains a PSC hypersurface $\Sigma$ which separates the faces of the band. This is due to the classical Schoen-Yau conformal descent argument adapted to $\mu$-bubbles as outlined by Gromov \cite[Section 5]{gromovfour}. We emphasise that the largeness of the band is not needed for existence of a separating $\mu$-bubble, however it is required for the $\mu$-bubble to be PSC. This argument works for all $3 \leq n \leq 6$. We present a proof for completeness in Section 2. 

The aim of the second step is to prove that if a $5$-dimensional band $M^4 \times [-1, 1]$ over a simply connected closed smooth $4$-manifold contains a separating PSC hypersurface $\Sigma \subset M \times [-1, 1]$, then $M$ admits an exotic copy which is PSC. We manage this following the same strategy as R\"{a}de \cite[Proposition 6.4]{rade23} (see also, \cite[Proposition 3.1]{schzen}) in the high dimensional setting, and carefully implementing it in dimension $4$. This involves considering the cobordism $(W; \partial_- W, \partial_+ W)$ between $\partial_- W = M$ and $\partial_+ W = \Sigma$ given by one of the connected components of $(M \times [-1, 1]) \setminus \Sigma$, and improving it by Wall's normal surgery techniques \cite{wallbook} so that the inclusion $M \hookrightarrow W$ induces an isomorphism on $\pi_1$.

In the higher dimensional case ($\dim(W) \geq 6$), one can also improve it to an isomorphism on $\pi_2$.  This would imply one can construct a handle decomposition of $(W, \partial_-W = M)$ containing no handles of index less than $3$. In that case, reversing the handle decomposition, one obtains $M$ from $\Sigma$ by surgeries of codimension $\geq 3$. In this \textit{ideal} scenario, $M$ is PSC if $\Sigma$ is PSC, due to the Gromov-Lawson-Schoen-Yau surgery theorem which states that PSC is preserved under surgeries of codimension $\geq 3$.

However, in the $\dim(W) = 5$, one can not guarantee such an isomorphism on $\pi_2$. Nevertheless, the inclusion $i:M \hookrightarrow W$ induces an injection $i_*: \pi_2(M)\rightarrow \pi_2(W)$ and the cokernel of $i_*$ is a $\ZZ[\pi_1]$-module generated by framed embedded $2$-spheres in $V$. Attaching these as trivial $2$-handles to $M$, we get a sub-cobordism between $\Sigma$ and $M\#^{k}(S^2 \times S^2)$ with only $1$- and $2$-handles, i.e., $M\#^{k}(S^2 \times S^2)$ is obtained from $\Sigma$ by only $0$- and $1$- surgeries. Since these are codimension $\geq 3$ surgeries performed on a $4$-manifold, this shows $M\#^{k}(S^2 \times S^2)$ is PSC for some $k$. In case when $M$ is simply connected we use the extension of the Gromov-Lawson-Stolz classification to dimension $4$ that becomes available upon considering PSC upto homeomorphism (Observation \ref{theoremCD}).  

Together, these give the width inequality. To see this, consider the band $M \times [-1,1]$ over a simply connected $4$-manifold $M$. By the $\mu$-bubble descent argument in Section 2, it contains a separating $\mu$-bubble $\Sigma$. If the band $M \times [-1, 1]$ is long enough, then $\Sigma$ is also PSC. Then, by the above, $M$ is PSC up to a homeomorphism. Therefore, if $M$ is not PSC upto homeomorphism, then the band $M \times [-1, 1]$ must satisfy the width inequality.

\subsection{Outline of paper} In Section 2 we recall definitions and review the existence and second variation of $\mu$-bubbles. We also provide here the Schoen-Yau descent argument adapted to $\mu$-bubbles. In Section 3, we prove several lemmas which culminate in the proof of Proposition \ref{prop-sephyppsc}, which is the main surgery result. In Section 4, we provide the proofs for Theorems \ref{theoremB}, \ref{theoremA}, and \ref{theoremC} as well as Observation \ref{theoremCD}.

\subsection*{Acknowledgments}  
The authors would like to thank Mike Miller Eismeier for comments on an early draft. They would also like to thank Jonathan Rosenberg for his interest in this work and for pointing out some inaccuracies in the previous draft. The second author would like to thank his advisor Mahan Mj for various fruitful discussions as well as Thorger Gei{\ss} for discussions on normal maps. The second author is supported by the Department of Atomic Energy, Government of India, under project no.12-R\&D-TFR-5.01-0500.

\section{Separating $\mu$-bubbles and width inequality}

In this section we introduce some terminology that will be used in the rest of the paper. We also review the $\mu$-bubble technique and collect some of its main results. These results are primarily from \cite{gromovfour,zhutams}. We have provided them for completeness.

 \subsection{Background} 
 The $\mu$-bubble method was introduced by Gromov in \cite[$5\frac{5}{6}$]{gromov96}. He explained its utility for problems in scalar curvature more recently in his \textit{Four Lectures} \cite[Section 5]{gromovfour}. He observed that one can modify the area functional, by adding an appropriate term, in such a way that the second variation of the modified functional is still amenable to the Schoen-Yau conformal descent argument. The minimizers of this modified functional are called $\mu$-bubbles and the technique is called the $\mu$-bubble method.  One key benefit of giving up minimality is that it is easier to show that $\mu$-bubbles separate the ambient space, i.e., in the minimization process they do not escape to infinity or to the boundary. This feature will become clear in the proof below. It is one of the reasons behind the success of the $\mu$-bubble method in several recent advances in the study of manifolds with positive scalar curvature (cf. \cite{clannals,cll,rade23}). 

\subsection{Preliminary definitions and properties}

\begin{definition}[Bands, faces, and width]
Let $M$ be a Riemannian manifold. A band over $M$ is the manifold $X = M \times [-1,1]$. The faces of the band are $\plusm = M \times \{+1\}$ and $\minm = M \times \{-1\}$. The width of the band is the distance between the faces, i.e., 
$$\mathrm{width}(X, g) = \dist_g\{\minm,\plusm\}$$
\end{definition}

\begin{definition}[Separating hypersurface]
A separating hypersurface in a band $X=M \times [-1,1]$ is an embedded hypersurface $\Sigma \subset X^\circ$ contained in the interior $X^\circ \subset X$ of the band such that there is no curve in $X \setminus \Sigma$ connecting the faces, i.e., if  $\gamma:[0,1]\to X$ is such that $\gamma(0) \in \minm$ and $\gamma(1) \in \plusm$, then $\gamma$ intersects $\Sigma$.
\end{definition}

We now formally define $\mu$-bubbles on a Riemannian band $X^{n+1} = M^n \times [-1, 1]$. 

\begin{definition}[$\mu$-bubbles]
Consider a function $h \in C^1(X^\circ)$ with the boundary conditions that $h$ goes to $\pm \infty$ on the faces $\partial_{\pm} X$, respectively. Fix a Cacciopoli set (i.e., a set of finite perimeter) $U_0$ with smooth boundary $\partial U_0 \subset X^\circ$ and $\minm \subset U_0$. For any Cacciopoli set $U$ such that $U \Delta U_0 \Subset X^\circ$, consider the functional
\begin{equation}\label{eq: mububblefunctional}
            \mathcal{A}_h(U) :=\mathcal{H}^{n}(\partial^*  U) - \int_{X} (\chi_{U}-\chi_{U_0}) h \,d\mathcal{H}^{n+1}
\end{equation}
Here, $\partial^*U$ denotes the reduced boundary of $U$ \cite[p. 42]{giusti}. If $\widetilde{U}$ is a minimizer of $\mathcal{A}_h$ with boundary $ \partial \widetilde{U} = \Sigma$, then we will call $\Sigma$ a $\mu$-bubble.
\end{definition}

Note that the case $h\equiv 0$, corresponding to minimal surfaces, is ruled out by the boundary conditions imposed on $h$, i.e., that it goes to $\pm \infty$ on the faces $\partial_{\pm} X$. Heuristically, this condition pushes the bubble away from the faces, giving a separating hypersurface. As explained earlier this is a key benefit of using $\mu$-bubbles, as a separating stable minimal hypersurface may not always exist.

The following was outlined by Gromov in \cite[Section 5]{gromovfour} and rigorously obtained by Zhu \cite{zhutams}.    

\begin{lemma}[Existence and Second Variation of $\mu$-bubbles]\label{lemm: secondvariation}
    Let $2\leq n \leq 6$ and $X^{n+1}$ be a Riemannian band. Then, $\mathcal{A}_h$ has a minimizer with smooth boundary  $\Sigma^{n}$ that separates the faces of the band. It has mean curvature $H_{\Sigma}=h|_{\Sigma}$ and the second variation formula is given by, 
    \begin{equation}
          \frac{1}{2}\int_{\Sigma} \left( R_M + |A|^2+ h^2 + 2\langle \nabla h,\nu\rangle \right) f^2 \,dvol_{\Sigma} \leq \int_{\Sigma} \left( |\nabla f|^2 + \frac{R_{\Sigma}}{2}f^2 \right)  \dvols 
    \end{equation}
    for all $f \in C^1_c(\Sigma)$.
\end{lemma}

Before we proceed further, we collect a small lemma. This result is behind the Schoen-Yau conformal descent.

\begin{lemma}\label{lem:kwlemma}
    Let $n\geq 3$ and $\Sigma^n$ be a closed Riemannian manifold. If 
 the operator $-\Delta_{g} + R/2$ is positive, i.e., $\lambda_1(-\Delta_{g} + R/2)>0$, then the metric $g$ can be conformally deformed to give a metric on $\Sigma$ with positive scalar curvature.  
\end{lemma}
\begin{proof}
For $n\geq 3$, consider the conformal Laplacian operator,
$$L_{\Sigma}  := - \Delta_{g} + \frac{n-2}{n-1}\frac{R}{4}$$ 
Note that $1/2 > (n-2)/(4n-4)>0$ if $n\geq 3$. Therefore, we have the following eigenvalue comparison
$$\lambda_1(L_{\Sigma})>\bigg(\frac{n-2}{2n-2} \bigg) \cdot \lambda_1(-\Delta_{g} + R/2) $$
If $\lambda_1(-\Delta_{g} + R/2)>0$, then the first eigenvalue of the conformal Laplacian, $\lambda_1(L_{\Sigma})>0$. Therefore, $\Sigma$ admits a PSC metric by a result of Kazdan and Warner \cite{kw75}. See also \cite[p. 9]{sydescent}.
\end{proof}
\begin{remark}
    The lemma is also true when $n=2$. There, it follows directly. One takes the test function $f \equiv 1$ in the defintion of $\lambda_1$ and uses Gauss-Bonnet theorem to obtain that $\Sigma$ is diffeomorphic to $S^2$.
\end{remark}
\subsection{Descent argument} We can now present the conformal descent argument for $\mu$-bubbles as outlined in \cite[Section 5]{gromovfour} and based on \cite[p. 20, Theorem 6.11]{chodoshnotes}. For $h\equiv 0$, this is exactly the argument of Schoen and Yau \cite{sydescent}. 

\begin{proposition}\label{prop:sepmububble}
    For $2 \leq n \leq 6$, let $X = M^{n} \times [-1,1]$ be a band over a closed manifold $M^{n}$. If $X$ has scalar curvature $R_X \geq n(n+1)\kappa^2$ and is such that 
    $$\dist_g(\plusm,\minm) > \frac{2\pi}{n+1}\frac{1}{\kappa}$$
    then $X$ contains a separating hypersurface that admits a PSC metric.
\end{proposition}
\begin{remark}
    The computation below can simplified by dropping the $|A|^2$ term. However, that does not provide a sharp width bound. 
\end{remark}
\begin{proof}
Given a function $h$ meeting the boundary conditions, the existence of a minimizing $\mu$-bubble $\Sigma$ that is a separating hypersurface for $X$ is contained in Lemma \ref{lemm: secondvariation}.  The rest of the argument is devoted to showing that $\Sigma$ is PSC.

We first recall the second variation formula from Lemma \ref{lemm: secondvariation}. Observe that by Cauchy-Schwarz $h^2/n \leq |A|^2$ on $\Sigma$ since $H_{\Sigma}=h|_{\Sigma}$. Therefore we get the following for any test function $f\in C_c^1(\Sigma)$:
\begin{equation}\label{eq:2ndvarinlemma}
      \int_{\Sigma} \bigg( \frac{R_X}{2} + \frac{n+1}{n} h^2 + \langle \nabla h,\nu\rangle \bigg) f^2  \leq \int_{\Sigma} |\nabla f|^2 + \frac{R_\Sigma}{2} f^2  
\end{equation}
Note here that if $h \equiv 0$, then we are in the minimal surface case (assuming existence). Then, Lemma \ref{lem:kwlemma} is immediately applicable. Below, we see that it is still applicable when band width is large enough, as then $h$ can be chosen such that the terms on RHS are bounded below by some $\delta >0$.
 
Denote the width by $L = \dist(\plusm,\minm)$ and denote by $\rho : X \to \RR$, a smoothening of the distance function from $\minm$. It is such that $\rho \equiv 0$ on $\minm$, $\rho \equiv 1$ on $\plusm$ and $|\nabla \rho| \leq 1$. Then, we can take the function $h$ in the $\mu$-bubble energy functional $\mathcal{A}_h$ to be 
$$h(x) = \frac{2n}{n+1}\frac{\pi}{\kappa^2L} \tan\bigg(\frac{\pi}{2}\bigg[\frac{2\rho(x)}{L}-1\bigg]\bigg)$$
Note that $h$ goes to $\pm \infty$ on the faces of the band as required. Further we have,
$$ 
    \nabla h = \frac{2n}{n+1}\frac{\pi^2}{\kappa^2L^2} \sec^2\bigg(\frac{\pi}{2}\bigg[\frac{2\rho(x)}{L}-1\bigg]\bigg) \nabla \rho
$$
Since $|\nabla \rho| \leq 1$ and $\sec^2(x)=1+\tan^2(x)$, this gives
$$ 
        -|\nabla h| \geq -\frac{2n}{n+1}\frac{\pi^2}{\kappa^2L^2} - \frac{n+1}{2n}h^2
$$
Therefore, the terms involving $h$ are bounded from below by a constant:
\begin{equation*}
    \frac{n+1}{2n}h^2 + \langle \nabla h,\nu\rangle \geq \frac{n+1}{2n}h^2 - |\nabla h| \geq -\frac{2n}{n+1}\frac{\pi^2}{\kappa^2L^2}
\end{equation*}
Plugging this back into the second variation formula \eqref{eq:2ndvarinlemma}, along with $R_X \geq n(n+1)\kappa^2$, we obtain
$$  \bigg( \frac{n(n+1)}{2} -\frac{2n}{n+1}\frac{\pi^2}{\kappa^2L^2} \bigg)\int_{\Sigma} f^2  \leq \int_{\Sigma} |\nabla f|^2 + (R_{\Sigma}/2)f^2 $$
Observe that if  $L= 2\pi/(n+1)\kappa$, then the left hand side is exactly $0$. Therefore, as per hypothesis, if the width $ L> 2\pi/(n+1)\kappa$, then the left hand side is bounded from below by some $\delta>0$. This gives,
$$0 < \delta \int_{\Sigma} f^2 <\int_{\Sigma}|\nabla f|^2 + (R_{\Sigma}/2)f^2$$
Therefore,
$$ 0 < \delta < \frac{\int_{\Sigma} |\nabla f|^2 + (R_{\Sigma}/2)f^2}{\int_{\Sigma}f^2}$$
Since this is true for all $f \in C^1_c(\Sigma)$, by taking an infimum over them we obtain:
$$0 < \delta <  \lambda_1(-\Delta_g + R/2) $$
Therefore, the $\mu$-bubble $\Sigma$ is PSC by Lemma \ref{lem:kwlemma}.  
\end{proof}

\section{Surgery arguments in dimension $4$}

The main goal of this section is to prove Proposition \ref{prop-sephyppsc}. To do that we first need to prove a few prepatory statements. The reader may wish to directly go to the proof of Proposition \ref{prop-sephyppsc}  and refer to the proofs of the lemmas as and when they are used therein. 

\subsection{Preliminaries on cobordism and handle decompositions}

We begin by fixing some notation that will be frequently used throughout this section.

\begin{definition}
    A \emph{cobordism} is a triple of smooth manifolds $(W; \partial_{-} W, \partial_{+} W)$ such that 
    $$\partial W = \partial_{-} W \sqcup \partial_{+} W$$
    $\partial_{-} W$ (resp. $\partial_{+} W$) is called the \emph{negative} (resp. \emph{positive}) boundary of $W$. $W$ is said to be a cobordism \emph{from} $\partial_{-} W$ \emph{to} $\partial_{+} W$. We shall often suppress this sense of direction, and simply say $W$ is a cobordism \emph{between} $\partial_{-} W$ and $\partial_{+} W$.
\end{definition}
 
Henceforth, all cobordisms will be assumed to be connected and between connected manifolds. Recall that an $n$-dimensional $k$-handle is $h^k := D^k \times D^{n-k}$. These are \emph{attached} to manifolds with boundaries $(X, \partial X)$ by diffeomorphisms $\phi : \partial D^k \times D^{n-k} \to \partial X$. The disk $D^k \times \{0\}$ (resp. $\{0\} \times D^{n-k}$) is called the \emph{core} (resp. \emph{co-core}) of $h^k$. The boundary of the core (resp. co-core) is called the \emph{attaching (resp. belt) sphere} of $h^k$.

\begin{definition}
    Given an $n$-dimensional cobordism $(W^n; \partial_{-} W, \partial_{+} W)$, a \emph{handle decomposition} of $W$ is a handle decomposition of $W$ rel. $\partial_{-} W$, i.e., a description of $W$ as obtained from $\partial_{-} W \times [0, 1]$ by attaching handles $D^k \times D^{n-k}$ of various indices $1 \leq k \leq n$ to $\partial_{-} W \times \{1\} \subset \partial_{-} W \times [0, 1]$. 

    Any such handle decomposition gives a dual handle decomposition of $W$ rel. $\partial_+ W$ in a canonical way, by reversing the roles of the core and co-core of each handle. We shall say this is obtained from \emph{reversing} the original handle decomposition of $W$. 
\end{definition}

The following lemma is well known, and is central to the proof of the $h$-cobordism theorem. We provide a proof for completeness.

\begin{lemma}[Arranging algebraically cancelling handles]\label{lem-algcancel} Let $(W^n; \partial_{-} W, \partial_+ W)$ be an $n$-dimensional cobordism, and let $k \leq n - 2$. Suppose 
\begin{enumerate}
\item $W$ admits a handle decomposition without handles of index $\leq k-1$. 
\item Moreover, $H_k(W, \partial_{-} W) = 0$
\end{enumerate}
Then, there exists a new handle decomposition of $W$ such that for every $k$-handle $h^k_i$, there is a $(k+1)$-handle $h^{k+1}_i$ so that the homological intersection number between the belt sphere of $h^k_i$ and the attaching sphere of $h^{k+1}_i$ is $+1$. \end{lemma}

\begin{proof}We use the language of handle homology, which is a variant of cellular homology for handle decompositions (see \cite[pg.~111]{gsbook}). As $H_k(W, \partial_{-} W) = 0$, every $k$-cycle in the handle chain complex of $(W, \partial_- W)$ bounds a $(k+1)$-chain. Therefore, for every $k$-handle $h^k_i$ of $(W, \partial_{-} W)$, there exists integers $c_{i, j} \in \ZZ$ such that
\begin{equation}\label{eq-1}[h^k_i] = \sum_j c_{i, j} \, d[h^{k+1}_j],\end{equation}
in the handle chain complex of $(W, \partial_- W)$. Here, the sum on the right hand side of Equation (\ref{eq-1}) runs over all $(k+1)$-handles of $W$. We modify the handle decomposition of $W$ as follows:
\begin{enumerate}
\item For every $k$-handle $h^k_i$ of $W$, create a canceling pair consisting of a $(k+1)$-handle and a $(k+2)$-handle. Let us denote these as $h^{k+1}_{i, 0}$ and $h^{k+2}_{i, 0}$, respectively.
\item Handle-slide $h^{k+1}_{i, 0}$ over each $(k+1)$-handle $h^{k+1}_j$, $c_{i, j}$ times. 
\end{enumerate}
Therefore, by construction,
$$[h^{k+1}_{i, 0}] = \sum_j c_{i, j} [h^{k+1}_j]$$
Consequently, $[h_i^k] = d[h^{k+1}_{i, 0}]$. Thus, in this new handle decomposition of $W$, the homological intersection number between the belt sphere of $h^k_i$ and the attaching sphere of $h^{k+1}_{i, 0}$ is $+1$, for all $i$. This finishes the proof.
\end{proof}

\subsection{Trading handles in cobordisms of dimension $5$}

Handle-trading is a standard technique used to simplify handle decompositions of cobordisms. It is well known that some of the techniques also extends to $5$-dimensional cobordisms with some effort, see \cite[Section 20.1]{detbook}. In this section we record some of the results relevant for our purpose. 

\begin{lemma}[Unknotting in dimension $4$]\label{lem-unknot} Let $M^4$ be a $4$-dimensional manifold and $\gamma \subset M$ be an embedded, nullhomotopic loop. Then $\gamma$ bounds a smoothly embedded disk in $M$.
\end{lemma}

\begin{proof}Since $\gamma$ is nullhomotopic, it must bound a clean immersed disk in $M$. Let $p_1, \cdots, p_n$ be the points of self-intersection. Join each $p_i$ to the boundary of the disk by an arc $\alpha_i$, $1 \leq i \leq n$. We apply the Whitney finger move along the arc $p_i$ to push the self-intersections out of the boundary of the disk (see \cite[Section 11.2]{detbook}, also \cite[Figure 11.4]{detbook}).
\end{proof}

\begin{lemma}[Trading $1$-handles for $3$-handles]\label{lem-trade} Let $(W; \partial_- W, \partial_+ W)$ be a $5$-dimensional cobordism, such that the inclusion $\partial_- W \hookrightarrow W$ is an isomorphism in $\pi_1$. Then, there exists a handle decomposition of $W$ without handles of index $0$ or $1$. 
\end{lemma}

\begin{proof}
Choose a handle decomposition of $W$. By connectedness, we may assume there are no $0$-handles. Let $V \subset W$ be the union of all $1$- and $2$-handles of $W$. We consider $V$ as a cobordism with negative boundary $\partial_{-} V = \partial_{-} W$, and positive boundary $\partial_+ V = \partial V \setminus \partial_- V$. Since $W$ is obtained by attaching handles of index $\geq 3$ to $V$, the inclusion $V \hookrightarrow W$ induces an isomorphism on $\pi_1$. Notice that the composition $\partial_- W = \partial_{-} V \hookrightarrow V \hookrightarrow W$
induces an isomorphism on $\pi_1$ by hypothesis. Therefore, the inclusion $\partial_- W \hookrightarrow V$ must also induce an isomorphism on $\pi_1$. Reversing the handle decomposition of $(V; \partial_{-} V, \partial_{+} V)$, we see $\partial_{+} V$ is obtained from $\partial_{-} V$ by attaching handles of index $\geq 3$. Therefore, $\partial_{+} V \hookrightarrow V$ induces an isomorphism on $\pi_1$, as well.

Let $h^1 \subset V$ be a $1$-handle. Let $\alpha \subset \partial_{+} V$ be an arc running parallel to the core of the $1$-handle along the boundary of $h^1$. The endpoints of $\alpha$ lie on the feet of the $1$-handle, which may be connected by an arc $\beta$ lying in $\partial_- W$. Next, since the inclusion
$$\partial_- W \hookrightarrow V$$
induces an isomorphism on $\pi_1$, the loop $[\alpha \cup \beta] \in \pi_1(V)$ must be the image of a loop $[\gamma] \subset \pi_1(\partial_- W)$. Let us denote $\beta' := \beta \# \gamma^{-1}$, given by interior connect sum of the arc $\beta$ with the loop $\gamma \subset \partial_{-} W$. Then, $\alpha \cup \beta' \subset \partial_{+} V \subset V$ must be nullhomotopic within $V$. However, from the previous paragraph, we know that the inclusion $$\partial_{+}V \hookrightarrow V,$$ also induces an isomorphism on $\pi_1$. Thus, $\alpha \cup \beta' \subset \partial_+ V$ is nullhomotopic in $\partial_{+} V$. 

By Lemma \ref{lem-unknot}, $\alpha \cup \beta' \subset \partial_{+} V$ bounds an embedded disk $D^2 \subset \partial_+ V$. We insert a layer $\partial_+ V \times I$ into the cobordism $W$, where $\partial_+ V \times I$ contains a $2$- and $3$-handle cancelling pair given by $h^2 = \nu(\partial(D^2 \times I) \setminus D^2 \times \{0\})$ and $h^3 = \nu(D^2 \times I)$. By construction, $h^2$ has attaching sphere $\alpha \cup \beta'$ which intersects the belt sphere of $h^1$ exactly once, geometrically. Therefore, $h^1$ and $h^2$ are a cancelling pair of handles. By cancelling them, we reduce the number of $1$-handles in $W$ by one. Proceeding inductively, we can likewise trade all the $1$-handles in $W$ for $3$-handles, as desired.
\end{proof}

\subsection{Surgery for normal maps in dimension $4$}

In this subsection, we extend a result due to R{\"a}de \cite[Proposition 6.4]{rade23} (see also \cite[Proposition 3.1]{schzen}) on making certain high dimensional cobordisms sufficiently homotopically connected, to cobordisms between manifolds of dimension $4$. This involves low dimensional implementation of standard ideas from high dimensional surgery theory of normal maps due to Wall \cite[Section 0.1]{wallbook}. We begin with the definition of normal maps.

\begin{definition}[Stable normal bundle and normal maps]
        The stable normal bundle of a smooth manifold $M$, denoted as $\nu(M)$, is the stable equivalence class\footnote{Two vector bundles $E, F$ over a space $X$ are said to be \emph{stably equivalent} if there exists $k, l \geq 0$ such that $E \oplus \varepsilon^k \cong F \oplus \varepsilon^l$. Here, $\varepsilon^k$ denotes the trivial vector bundle of rank $k$.} of the normal bundle of a Whitney embedding of $M$ in $\RR^n$, for some $n \geq 1$. 

        Let $M, N$ be smooth manifolds. A map $f : M \to N$ is said to be \emph{normal} if $f^*\nu(N) \oplus TM$ is stably trivial.
\end{definition} 

\begin{lemma}[Making $5$-dimensional cobordisms $1$-connected]\label{lem-2conn} Let $(W; \partial_{-} W, \partial_+ W)$ be a $5$-dimensional cobordism, and let $r : W \to \partial_{-} W$ be a retract which is a normal map. There exists a new $5$-dimensional cobordism $V$ with $\partial_{\pm} V = \partial_{\pm} W$ such that \begin{enumerate}
    \item The inclusion $i : \partial_{-} V \hookrightarrow V$ induces an isomorphism on $\pi_1$ and an injection on $\pi_2$.
    \item Any embedded $2$-sphere in the interior of $V$ representing a class in the cokernel of $i_* : \pi_2(\partial_- V) \to \pi_2(V)$ has trivial normal bundle.
\end{enumerate} 
\end{lemma}

\begin{proof}
Since $r$ is a retract, $r_* : \pi_1(W) \to \pi_1(\partial_{-} W)$ is a surjection. As $W, \partial_{-} W$ are compact, the groups $\pi_1(W)$ and $\pi_1(\partial_{-} W)$ are finitely presented. Therefore, $\ker(r_*) < \pi_1(W)$ is \emph{normally} finitely generated. We choose a generating set for $\ker(r_*)$, so that:
$$\ker(r_*) = \langle \! \langle [\gamma_1], \cdots, [\gamma_n] \rangle\! \rangle$$
Since $\dim(W) = 5 > 2$, we may ensure by a homotopy that the representatives $\gamma_i \subset W^\circ$ are embedded loops strictly contained in the interior of $W$. 

Let $N(\gamma_i) \subset W^\circ$ be a normal neighborhood. Since $r : W \to \partial_-W$ is a normal map, $r^* \nu(\partial_- W) \oplus TW$ is stably trivial. Therefore, the restriction of $r^* \nu(\partial_- W) \oplus TW$ to $\gamma_i \subset W^\circ$ is also stably trivial. As $r|_{\gamma_i}$ is null-homotopic, there exists a map $f : D^2 \to \partial_- W$ such that $f|_{\partial D^2} = r|_{\gamma_i}$. Thus,
$$r^*\nu(\partial_- W)|_{\gamma_i} = g^*\nu(\partial_- W)|_{\partial D^2}$$
However, $g^*\nu(\partial_-W)$ is trivial as it is a bundle over $D^2$. Thus, $r^*\nu(\partial_- W)|_{\gamma_i}$ is trivial as well. Consequently, $TW|_{\gamma_i}$ must be stably trivial. However,
$$TW|_{\gamma_i} = T\gamma_i \oplus N(\gamma_i) \cong TS^1 \oplus N(\gamma_i)$$
As $TS^1 \cong S^1 \times \RR$ is trivial, $N(\gamma_i)$ must be stably trivial. If $N(\gamma_i)$ is not orientable, then the first Stiefel-Whitney class must be non-zero. Since that would contradict stable triviality, we conclude it is orientable. As oriented vector bundles over the circle are trivial, $N(\gamma_i)$ is trivial. Hence, $N(\gamma_i) \cong S^1 \times D^4$.

We perform a surgery along these loops in the interior of $W$, by cutting out the normal neighborhoods $N(\gamma_i) \cong S^1 \times D^4$ and gluing back $D^2 \times S^3$. Let us call the resulting manifold $V$, which is a new cobordism between $\partial_\pm V = \partial_\pm W$. Since $r|_{\gamma_i}$ is a nullhomotopic, $r$ extends over the surgery to a retract $r' : V \to \partial_{-} V$. By construction, $(r')_* : \pi_1(V) \to \pi_1(\partial_- V)$ is an isomorphism. As $r'$ is also a retract, $(r')_* : \pi_2(V) \to \pi_2(\partial_- V)$ is a surjection. Since the inclusion map $i : \partial_- V \hookrightarrow V$ satisfies $r \circ i = \mathrm{id}_{\partial_- V}$, $i$ must induce an isomorphism on $\pi_1$ and an injection on $\pi_2$. This proves Part $(1)$.

Next, we proceed to prove Part $(2)$. Note that as $r'$ is a retract, there is a splitting:
$$\pi_2(V) \cong \ker((r')_*) \oplus \mathrm{im}(i_*)$$
Therefore, $\mathrm{coker}(i_*) \cong \ker((r')_*)$. Let $\sigma \subset V^\circ$ be an embedded $2$-sphere in the interior of $V$ such that $[\sigma]$ is in the kernel of $(r')_* : \pi_2(V) \to \pi_2(\partial_- V)$. Since $r : W \to \partial_{-} W$ is a normal map and $r' : V \to \partial_- V$ is obtained from $1$-surgeries to $r$, $r'$ is also a normal map\footnote{Performing surgery below the middle dimension on a normal map also produces a normal map. For details, see \cite[Section 0.1]{wallbook}}. Thus, $(r')^*\nu(\partial_- V) \oplus TV$ is stably trivial. Therefore, the restriction of $(r')^*\nu(\partial_- V) \oplus TV$ to $\sigma \subset V$ is also stably trivial. Since $r'|_{\sigma}$ is nullhomotopic, there exists a map $g : D^3 \to \partial_- V$ such that $g|_{\partial D^3} = r'|_{\sigma}$. Thus,
$$(r')^* \nu(\partial_- V)|_{\sigma} = g^* \nu(\partial_{-} V) |_{\partial D^3}$$
However, $g^* \nu(\partial_{-} V)$ is trivial as it is a bundle over $D^3$. Thus, $(r')^* \nu(\partial_- V)$ is trivial. Consequently, $TV|_{\sigma}$ must be stably trivial. However, 
$$TV|_{\sigma} = T\sigma \oplus N(\sigma) \cong TS^2 \oplus N(\sigma).$$
Since $TS^2$ is stably trivial, this in turn forces $N(\sigma)$ to be stably trivial. Note that $N(\sigma)$ is a rank $3$ bundle over $\sigma \cong S^2$. Since $BO(3) \hookrightarrow BO(\infty)$ induces an isomorphism on $\pi_2$, $N(\sigma)$ must be in fact be trivial, i.e. $N(\sigma) \cong S^2 \times D^3$.
\end{proof}

\subsection{Obtaining a PSC metric on a stabilization}

\begin{proposition}\label{prop-sephyppsc} Let $M$ be a closed smooth $4$-manifold. Let $\Sigma \subset M \times [-1, 1]$ be an embedded hypersurface separating $M \times \{-1\}$ and $M \times \{1\}$. If $\Sigma$ admits a positive scalar curvature metric, then for some $k \geq 0$, $M \#^k (S^2 \times S^2)$ admits a positive scalar curvature metric.
\end{proposition}

\begin{proof}
Let $W \subset M \times [-1, 1]$ denote the compact domain cobounded by $M \times \{-1\}$ and $\Sigma$. We think of $W$ as a cobordism from $\partial_{-} W = M$ to $\partial_{+} W = \Sigma$. Restricting the projection $M \times [-1, 1] \to M \times \{-1\}$ to $W$, we conclude $W$ admits a retract $r : W \to \partial_{-} W = M$. Furthermore, since the projection $M \times [-1, 1] \to M$ is a normal map, so is the restriction $r : W \to M$. Using Lemma \ref{lem-2conn}, we find a new cobordism $(V; M, \Sigma)$ such that the inclusion $i : M \hookrightarrow V$ induces an isomorphism on $\pi_1$ and an injection on $\pi_2$, and all embedded $2$-spheres in the interior of $V$ representing a class in the cokernel of $i_* : \pi_2(M) \to \pi_2(W')$ has trivial normal bundle in $V$.

Since $i$ induces an isomorphism on $\pi_1$, we may choose a handle decomposition of $(V; M, \Sigma)$ without $1$-handles, by Lemma \ref{lem-trade}. Let $h^2_1, \cdots, h^2_k$ be all the $2$-handles in $V$. Let $\delta > -1$ be a constant. We suppose $h^2_i$ are attached to $M \times [-1, \delta] \subset V$ along the level set $M \times \{\delta\}$. Using the fact that $i$ induces an isomorphism on $\pi_1$ once again, we observe that the attaching circles of $h^2_i$ must be null-homotopic in $M \times \{\delta\}$. We choose disjoint embedded disks bounding these null-homotopic circles within the $5$-manifold $M \times [-1, \delta]$. By capping off the core of the $2$-handles $h^2_i$ by these disks, we obtain disjoint embedded spheres $S^2_i \subset V$. Let us call $S^2_i$ the \emph{core sphere} of the $2$-handle $h^2_i$.

Let us denote $\pi := \pi_1(M) \cong \pi_1(V)$. Observe that the $2$-skeleton $V^{(2)}$ of $V$ is homotopy equivalent to $M \vee (\vee_{i = 1}^k S^2_i)$. Let $z_i = [S^2_i]$ denote the homotopy classes of the core spheres. Then, we have an isomorphism of $\ZZ[\pi]$-modules:
$$\pi_2(V^{(2)}) \cong H_2(\widetilde{M \vee (\vee_{i = 1}^k S^2_i)}) \cong \pi_2(M) \oplus \left ( \bigoplus_{i = 1}^k \ZZ[\pi]\langle z_i \rangle \right )$$
The map on $\pi_2$ induced by the inclusion $M \hookrightarrow V^{(2)}$ is an isomorphism onto the first factor above. Since $\pi_2(V) \cong \pi_2(V^{(3)})$ and $V^{(3)}$ is obtained from $V^{(2)}$ by attaching $3$-cells, there exists a submodule $A \subset \pi_2(V^{(2)})$ such that $\pi_2(V) \cong \pi_2(V^{(2)})/A$. Note that $A$ must intersect $\pi_2(M) \subset \pi_2(V^{(2)})$ trivially. Otherwise, a non-zero element in the intersection will live in the kernel of the composition $i : M \hookrightarrow V^{(2)} \hookrightarrow V$. This is absurd, as $i$ induces an injection on $\pi_2$. Therefore, $A \subset \oplus_{i = 1}^k \ZZ[\pi]\langle z_i \rangle$. Consequently, 
$$\pi_2(V) \cong \pi_2(M) \oplus \left (\bigoplus_{i = 1}^k \ZZ[\pi]\langle z_i \rangle \right )/A$$
Therefore, the cokernel of $i_*$ is generated by $z_1, \cdots, z_k$. Hence, all the core spheres $S^2_i \subset V$ of the $2$-handles of $V$ have trivial normal bundle, by Part $(2)$ of Lemma \ref{lem-2conn}.

We pass to the level set in $V$ after attaching all the $2$-handles $h^2_1, \cdots, h^2_k$ to $M$. Since the $2$-handles are attached along nullhomotopic curves and their core spheres are framed, the result is diffeomorphic to $M \#^k (S^2 \times S^2)$. The remaining sub-cobordism $V' \subset V$ between $M \#^k (S^2 \times S^2)$ and $\Sigma$ consists of only $3$- and $4$-handles. Reversing the cobordism $V'$, we see that $M \#^k (S^2 \times S^2)$ is obtained by attaching $1$- and $2$-handles to $\Sigma$. In other words, $M \#^k (S^2 \times S^2)$ is obtained from $0$- and $1$-surgeries on $\Sigma$. Since these are codimension $\geq 3$ surgeries performed on a PSC manifold $\Sigma$, the resulting manifold must continue to be PSC by \cite{glsc, sydescent}. Thus, $M \#^k (S^2 \times S^2)$ is PSC. 
\end{proof}

\section{Width inequality and $S^1$-stability in dimension $4$ }

\subsection{General case} 
We first prove a version of the width inequality and $S^1$-stability for all $4$-manifolds. 

\begin{theorem}\label{thm:stablypscband}
    Let $M$ be a closed smooth $4$-manifold such that $M\#^k(S^2 \times S^2)$ does not admit a PSC metric for any $k\geq 0$. Let $g$ be a metric on the band $X = M^4 \times [-1,1]$ such that $R_g \geq 20\kappa^2$. Then,
$$\wid(X^5,g) := \dist_g(M \times \{-1\},M \times \{1\}) \leq \frac{2\pi}{5\kappa}$$
 \end{theorem}
\begin{proof}
    First note that due to Proposition \ref{prop-sephyppsc}, the band $X$ does not admit any PSC hypersurface that separates the faces. Therefore, if the width inequality is not satisfied, then $X$ satisfies all hypothesis of Proposition \ref{prop:sepmububble} with $n=4$. Consequently, $X$ admits a separating $\mu$-bubble that is PSC. Contradiction.  
\end{proof}
\begin{corollary}
     Let $M^4$ be any closed smooth $4$-manifold. If $\mts$ is PSC then $M\#^k(S^2 \times S^2)$ is PSC for some $k\geq 0$.
\end{corollary}
\begin{proof}
     Since, $\mts$ is PSC, we may lift the metric to a complete and uniformly PSC metric on the cover $M \times \RR$. Therefore, there is a $\kappa>0$ such that $R_{M \times \RR}\geq 20\kappa^2$. Next, we take $t$ large enough so that the width of $M\times[-t,t]$ is greater than $2\pi/5\kappa$. By Theorem \ref{thm:stablypscband} this implies that $M\#^k(S^2 \times S^2)$ is PSC for some $k\geq 0$.
\end{proof}

\subsection{Simply connected case}
We first extend the Gromov-Lawson-Stolz classification of simply connected manifolds with PSC to dimension $4$.
\begin{observation}[Gromov-Lawson-Stolz in $n=4$]\label{obs-gls4d}
     Let $M^4$ be a closed simply connected smooth $4$-manifold. If $M$ is not spin, then it is PSC upto homeomorphism. If $M$ is spin, then it is PSC upto homeomorphism if and only if $\hat{A}(M)=0$.
\end{observation}
\begin{proof}
Suppose $M$ is not spin. Since $M$ is a simply connected this is equivalent to $M$ having with an odd intersection form $Q_M$. Therefore, with $m=b_2^+(M)$ and $n=b_2^-(M)$, we have that $Q_M=m\langle 1 \rangle \oplus n\langle -1 \rangle $. Consequently, $M$ is homeomorphic to $m\cptwo \# n\cptwobar$ by Freedman's theorem \cite{4dpc}. The latter is a PSC manifold.

Suppose $M$ is spin. By the classical result of Lichnerowicz \cite{lic63}, if $M$ is spin and $\hat{A}(M) \neq 0$, then $M$ does not admit a PSC metric up to homeomorphism. This is true upto homeomorphism because $\hat{A}$ is a homotopy invariant (in particular, a homeomorphism invariant) in dimension $4$. In the other direction, it is noted earlier that in dimension $4$, $\hat{A}(M)=-\sigma(M)/8$, where $\sigma=b_2^+-b_2^-$ is the signature. Therefore, if $\hat{A}(M)=0$, this implies that $M$ has zero signature. But all simply connected spin $4$-manifolds with zero signature are homeomorphic to a connect sum $\#^{k}(S^2\times S^2)$ for some $k\geq 1$ \cite{4dpc}. The latter is a PSC manifold. 
\end{proof}

\begin{proposition}\label{cor-sephyppsc}
    Let $M$ be a closed simply connected smooth $4$-manifold. Let $\Sigma \subset M \times [-1, 1]$ be an embedded hypersurface separating $M \times \{-1\}$ and $M \times \{1\}$. If $\Sigma$ admits a positive scalar curvature metric, then $M$ is PSC upto homeomorphism.
\end{proposition}

\begin{proof}
   Since, $M$ is simply connected by Observation \ref{obs-gls4d}, we need only consider the case where $M$ is spin. Then, $M \#^k (S^2 \times S^2)$ is a simply connected spin manifold, and by Proposition \ref{prop-sephyppsc}, it is PSC for some $k$. Therefore, $\widehat{A}(M \#^k (S^2 \times S^2)) = 0$. However,
    $$\widehat{A}(M \#^k (S^2 \times S^2)) = -\frac18 \sigma(M \#^k(S^2 \times S^2)) = -\frac18 \sigma(M) = \widehat{A}(M)$$
    Therefore, $\widehat{A}(M) = 0$. By Observation \ref{obs-gls4d}, $M$ is PSC upto homeomorphism, as desired.
\end{proof}

We can now establish the band width inequality for simply connected $4$-manifolds that are not PSC upto homeomorphisms. 

\begin{theorem}
    \label{thm:4dwidth}
    Let $M^4$ be a closed simply connected smooth $4$-manifold that is not PSC upto homeomorphism. Let $g$ be a metric on the band $X = M^4 \times [-1,1]$ such that $R_g \geq 20\kappa^2$. Then,
$$\wid(X^5,g) := \dist_g(M \times \{-1\},M \times \{1\}) \leq \frac{2\pi}{5\kappa}$$
\end{theorem}
\begin{proof}
  First note that by Proposition \ref{prop-sephyppsc}, the band $X$ over $M$ does not admit a PSC separating hypersurface. If the inequality is not satisfied, then the band $X$ satisfies all hypothesis of Proposition \ref{prop:sepmububble} with $n=4$. Consequently, $X$ admits a separating $\mu$-bubble that is PSC. Contradiction. 
\end{proof}

Finally, we prove $S^1$-stability for simply connected $4$-manifolds. 

\begin{theorem}Let $M^4$ be a simply connected closed smooth $4$-manifold. Then $M$ is PSC up to homeomorphism if and only if $M \times S^1$ is PSC. 
\end{theorem}
\begin{proof}
($\Longleftarrow$) First, assume that $\mts$ is PSC. Then we may lift the metric to a complete and uniformly PSC metric on the cover $M \times \RR$. Therefore, there is a $\kappa>0$ such that $R_{M \times \RR}\geq 20\kappa^2$. Next, we take $t$ large enough so that the width of $M\times[-t,t]$ is greater than $2\pi/5\kappa$. By Theorem \ref{thm:4dwidth}, this implies that M is PSC upto homeomorphism.

($\Longrightarrow$) Now, we assume that $M$ is PSC upto homeomorphism. So, there is a smooth $4$-manifold $M'$ which is PSC and $M$ is homeomorphic to $M'$. Note that, as $M'$ is PSC, so is $M'\times S^1$ (one can take the metric $g_{M'}+dt^2$). 
Now, since $M$ and $M'$ are simply connected, there is an $h$-cobordism $(W;M',M)$. Upon taking a product with $S^1$, we obtain an $h$-cobordism $(W\times S^1; M' \times S^1, \mts)$. Note that $\pi_1(\mts)=\ZZ$ as $M$ is simply connected. But the Whitehead group $\mathrm{Wh}(\ZZ)$ vanishes, therefore $\mts$ and $M' \times S^1$ are diffeomorphic by the $s$-cobordism theorem. Pulling back the PSC metric on $M' \times S^1$ by this diffeomorphism, we obtain $M \times S^1$ is also PSC.
\end{proof}

\bibliography{bib}
\bibliographystyle{alpha}
\end{document}